\documentclass[11pt,a4paper,fullpage]{article}
\usepackage[english]{babel}
\usepackage[utf8]{inputenc}
\usepackage{fancyhdr}
\usepackage[numbers]{natbib}

\usepackage{graphicx}
\usepackage{amsmath}
\usepackage{extsizes}
\usepackage{relsize}
\usepackage{multicol}
\usepackage{ragged2e}
\usepackage{fontenc}
\usepackage{mathcomp}
\usepackage{latexsym}
\usepackage{amssymb}
\usepackage{times}
\newtheorem{Teorema}{Theorem}[section]

\newtheorem{Posledica}[Teorema]{Corollary}

\newtheorem{Lema}[Teorema]{Lemma}
\newtheorem{Primedba}[Teorema]{Remark}
\usepackage{dsfont}
\newcommand\scalemath[2]{\scalebox{#1}{\mbox{\ensuremath{\displaystyle #2}}}} 
\numberwithin{equation}{section}
\hyphenchar\font=-1
\begin{document}
	\title {Fredholmness and Weylness of block operator matrices}	

\author{Nikola Sarajlija\footnote{corresponding author: Nikola Sarajlija; University of Novi Sad, Faculty of Sciences, Novi Sad 21000, Serbia; {\it e-mail}: {\tt nikola.sarajlija@dmi.uns.ac.rs}}\footnote{The author is supported by the Ministry of Education, Science and Technological Development of the Republic of Serbia under grant no.  451-03-68/2022-14/200125.}}
\maketitle

\begin{abstract}
This paper has aim to characterize Fredholmness and Weylness of upper triangular operator matrices having arbitrary dimension $n\geq 2$. We present various characterization results in the setting of infinite dimensional Hilbert spaces, thus extending some known results from Cao X. et al. (Acta Math. Sin. (Engl. Ser.) \textbf{22} (2006), no. 1, 169–178 and J. Math. Anal. Appl. \textbf{304} (2005), no. 2, 759–771)  and Zhang et al. (J. Math. Anal. Appl. \textbf{392} (2012), no. 2, 103–110) to the case of arbitrary dimension  $n\geq2$. We pose our results without using separability assumption, thus improving perturbation results from Wu X. et al. (Ann. Funct. Anal. \textbf{11} (2020), no. 3, 780–798 and Acta Math. Sin. (Engl. Ser.) 36 (2020), no. 7, 783–796).
\end{abstract}
\textit{$2020$ Math. Subj. Class:} 47A08, 47A53, 47A55, 47A05, 47A10.

\vspace{2mm}

\section{Introduction and notation}

This article is concerned with partial upper triangular operator matrices of arbitrary dimension $n\geq 2$. Term "partial" means that some of the entries of a matrix are given, while the others are unknown. These operators arise naturally in many research areas of operator theory. Indeed, suppose that an operator $T$ acting on a Hilbert space $\mathcal{H}$ is studied with respect to an orthogonal decomposition of $\mathcal{H}$. In other words, let ${M}$ be a closed subspace of $\mathcal{H}$. Then we have $\mathcal{H}=M\oplus M^\bot$ and $T=\begin{bmatrix}T_1 & T_2\\ T_3 & T_4\end{bmatrix},$  where $T_1:M\rightarrow M$,  $T_2:M^\bot\rightarrow M$,  $T_3:M\rightarrow M^\bot$,  $T_4:M^\bot\rightarrow M^\bot$. Among all closed subspaces of $\mathcal{H}$ we distinguish those which are invariant for $T$. If $M$ is a such subspace, then $T(M)\subseteq M$ and $T$ takes the upper triangular form $T=\begin{bmatrix}T_1 & T_2\\ 0 & T_4\end{bmatrix}$. We can implement this reasoning for an arbitrary number of closed subspaces, thus obtaining upper triangular operators of dimension $n>2$. Spectral properties of block operators are extensively studied by numerous authors (see \cite{CAO}, \cite{CAO2}, \cite{OPERATORTHEORY}-\cite{KOLUNDZIJA3}, and \cite{WU}-\cite{ZHANG}) and we continue this study for upper triangular operators of an arbitrary dimension.

Article is organized as follows. In the rest of this section we give notation and some basics on Fredholm theory with a few auxiliary results. Sections 2 and 3 deal with characterizing Weylness and Fredholmness in the setting of arbitrary infinite dimensional Hilbert spaces, respectively.

Let $\mathcal{H},\mathcal{H}_1,...,\mathcal{H}_n$ be Hilbert spaces. We use notation $\mathcal{B}(\mathcal{H}_i,\mathcal{H}_j)$ for the collection of all linear and bounded operators from $\mathcal{H}_i$ to $\mathcal{H}_j$, where $\mathcal{B}(\mathcal{H}):=\mathcal{B}(\mathcal{H},\mathcal{H})$ for short. If $T\in\mathcal{B}(\mathcal{H}_i,\mathcal{H}_j)$, by $\mathcal{N}(T)$ and $\mathcal{R}(T)$ we denote the kernel and the range of $T$, respectively. It is well known that $\mathcal{N}(T)$ is closed, thus complemented in Hilbert space $\mathcal{H}_i$. By $T^*\in\mathcal{B}(\mathcal{H}_j,\mathcal{H}_i)$ we denote the adjoint of $T$.

Let $D_1\in\mathcal{B}(\mathcal{H}_1),\ D_2\in\mathcal{B}(\mathcal{H}_2),...,D_n\in\mathcal{B}(\mathcal{H}_n)$ be given. Partial upper triangular operator matrix of dimension $n$ is
\begin{equation}\label{OSNOVNI}
T_n^d(A)=
\begin{bmatrix} 
    D_1 & A_{12} & A_{13} & ... & A_{1,n-1} & A_{1n}\\
    0 & D_2 & A_{23} & ... & A_{2,n-1} & A_{2n}\\
    0 &  0 & D_3 & ... & A_{3,n-1} & A_{3n}\\
    \vdots & \vdots & \vdots & \ddots & \vdots & \vdots\\
    0 & 0 & 0 & ... & D_{n-1} & A_{n-1,n}\\
    0 & 0 & 0 & ... & 0 & D_n      
\end{bmatrix}\in\mathcal{B}(\mathcal{H}_1\oplus \mathcal{H}_2\oplus\cdots\oplus \mathcal{H}_n),
\end{equation}
where $A:=(A_{12},\ A_{13},...,\ A_{ij},...,\ A_{n-1,n})$ is an operator tuple consisting of unknown variables $A_{ij}\in\mathcal{B}(\mathcal{H}_j,\mathcal{H}_i)$, $1\leq i<j\leq n,\ n\geq2$. We denote by $\mathcal{B}_n$ the collection of all described tuples $A=(A_{ij})$. This is a notation that originated in \cite{WU}, and was used later in \cite{WU2}, \cite{WU3}. One easily verifies that for $T_n^d(A)$ having form (\ref{OSNOVNI}), its adjoint operator matrix $T_n^d(A)^*$ is given by
\begin{equation}\label{ADJUNGOVANI}
T_n^d(A)^*=
\scalemath{0.9}{
\begin{bmatrix} 
    D_1^* & 0 & 0 & ... & 0 & 0\\
    A_{12}^* & D_2^* & 0 & ... & 0 & 0\\
    A_{13}^* &  A_{23}^* & D_3^* & ... & 0 & 0\\
    \vdots & \vdots & \vdots & \ddots & \vdots & \vdots\\
    A_{1,n-1}^* & A_{2,n-1}^* & A_{3,n-1}^* & ... & D_{n-1}^* & 0\\
    A_{1n}^* & A_{2n}^* & A_{3n}^* & ... & A_{n-1,n}^* & D_n^*      
\end{bmatrix}}\in\mathcal{B}(\mathcal{H}_1\oplus \mathcal{H}_2\oplus\cdots\oplus \mathcal{H}_n).
\end{equation}

Topic that interests us is a part of Fredholm theory. Here we give some basics, using notation from \cite{ZANA}. Let $T\in\mathcal{B}(\mathcal{H})$, $\alpha(T)=\dim\mathcal{N}(T)$ and $\beta(T)=\dim\mathcal{H}/\mathcal{R}(T)$. Quantities $\alpha(T)$ and $\beta(T)$ are called the nullity and the deficiency of $T$, respectively, and in the case where at least one of them is finite we define $\mathrm{ind}(T)=\alpha(T)-\beta(T)$ to be the index of $T$. Families of upper and lower semi-Fredholm operators, respectively, are defined as
$$
\begin{aligned}
\Phi_+(\mathcal{H})=\left\{ T\in\mathcal{B}(\mathcal{H}):\ \alpha(T)<\infty\ and\ \mathcal{R}(T)\ is \ closed\right\}
\end{aligned}
$$
and
$$
\begin{aligned}
\Phi_-(\mathcal{H})=\lbrace T\in\mathcal{B}(\mathcal{H}): \beta(T)<\infty\rbrace.
\end{aligned}
$$

The set of Fredholm operators is
$$\Phi(\mathcal{H})=\Phi_+(\mathcal{H})\cap\Phi_-(\mathcal{H})=\lbrace T\in\mathcal{B}(\mathcal{H}): \alpha(T)<\infty\ and\ \beta(T)<\infty\rbrace.$$
Put
$$\Phi_+^-(\mathcal{H})=\lbrace T\in\Phi_+(\mathcal{H}):\ \mathrm{ind}(T)\leq0\rbrace$$
and
$$\Phi_-^+(\mathcal{H})=\lbrace T\in\Phi_-(\mathcal{H}):\ \mathrm{ind}(T)\geq0\rbrace.$$
The previous two collections are called the sets of upper and lower Weyl operators, respectively.

Corresponding spectra of an operator $T\in\mathcal{B}(\mathcal{H})$ are defined as follows:\\
-the upper semi-Fredholm spectrum: $\sigma_{SF+}(T)=\lbrace\lambda\in\mathds{C}: \lambda-T\not\in\Phi_{+}(\mathcal{H})\rbrace$;\\
-the lower semi-Fredholm spectrum: $\sigma_{SF-}(T)=\lbrace\lambda\in\mathds{C}: \lambda-T\not\in\Phi_{-}(\mathcal{H})\rbrace$;\\
-the essential spectrum: $\sigma_{e}(T)=\lbrace\lambda\in\mathds{C}: \lambda-T\not\in\Phi(\mathcal{H})\rbrace$;\\
-the upper semi-Weyl spectrum: $\sigma_{aw}(T)=\lbrace\lambda\in\mathds{C}: \lambda-T\not\in\Phi_+^-(\mathcal{H})\rbrace$;\\
-the lower semi-Weyl spectrum: $\sigma_{sw}(T)=\lbrace\lambda\in\mathds{C}: \lambda-T\not\in\Phi_-^+(\mathcal{H})\rbrace$;\\[1mm]
We write $\rho_{SF+}(T), \rho_{SF-}(T), \rho_e(T), \rho_{aw}(T), \rho_{sw}(T)$ for their complements, respectively.

We list some elementary results from functional analysis.

\begin{Lema}\label{VEZA}
For $T\in\mathcal{B}(\mathcal{H})$ with closed range the following holds:\\
$(a)$ $\alpha(T)=\beta(T^*),\beta(T)=\alpha(T^*)$;\\
$(b)$ $T\in\Phi_+(\mathcal{H})$ if and only if $T^*\in\Phi_-(\mathcal{H})$;\\
$(c)$ $T\in\Phi_-(\mathcal{H})$ if and only if $T^*\in\Phi_+(\mathcal{H})$;\\
$(d)$ $\mathrm{ind}(T^*)=-\mathrm{ind}(T).$
\end{Lema}

\begin{Lema}\label{POMOCNALEMA}
Let $T_n^d(A)\in\mathcal{B}(\mathcal{H}_1\oplus\cdots\oplus\mathcal{H}_n).$ Then:
\begin{itemize}
\item[(i)] $\sigma_{SF+}(D_1)\subseteq\sigma_{SF+}(T_n^d(A))$;
\item[(ii)] $\sigma_{SF-}(D_n)\subseteq\sigma_{SF-}(T_n^d(A))$.
\end{itemize}
\end{Lema}

\begin{Lema}\label{DODATNALEMA}
Let $T\in\mathcal{B}(\mathcal{H})$. Then $\mathcal{R}(T)$ is closed if and only if $\mathcal{R}(T^*)$ is closed.
\end{Lema}

\begin{Teorema}(\cite{CARADUS})\label{KATOZATVOREN}
For $T\in\mathcal{B}(\mathcal{H})$ the following implication holds:
$$
\beta(T)<\infty\Rightarrow\mathcal{R}(T)\ is\ closed\ (thus\ complemented\ in\ \mathcal{H}).
$$
\end{Teorema}

\section{Weylness of $T_n^d(A)$}

In this section we assume $\mathcal{H}_1,...,\mathcal{H}_n$ to be arbitrary infinite dimensional Hilbert space. We generalize results of \cite{CAO2},\cite{ZHANG} from $n=2$ to an arbitrary dimension of upper triangular operators, and we pose perturbation results of \cite{WU3} without assuming separability of underlying spaces. 

\subsection{Weylness of $T_n^d(A)$}

We start with a result which deals with the upper Weyl spectrum of $T_n^d(A)$. 
\begin{Teorema}\label{LEVIVEJL'}
Let $D_1\in\mathcal{B}(\mathcal{H}_1),\ D_2\in\mathcal{B}(\mathcal{H}_2),...,D_n\in\mathcal{B}(\mathcal{H}_n)$. Consider the following statements:\\
$(i)$ $(a)$ $D_1\in\Phi_+(\mathcal{H}_1)$;\\
\hspace*{5.5mm}$(b)$ $\mathcal{R}(D_s)$ is closed for $2\leq s\leq n$ and 
\begin{equation}\label{PRVA}
\begin{aligned}
\Big(\alpha(D_s)\leq\beta(D_{s-1})\quad for\ 2\leq s\leq n,\\
\sum_{s=1}^n\beta(D_s)=\infty\Big)
\end{aligned}
\end{equation}
or $\Big(D_s\in\Phi_+(\mathcal{H}_s)$ for $2\leq s\leq n$ and $\sum\limits_{s=1}^n\alpha(D_s)\leq\sum\limits_{s=1}^n\beta(D_s)$\Big);\\[3mm]
$(ii)$ There exists $A\in\mathcal{B}_n$ such that $T_n^d(A)\in\Phi_+^-(\mathcal{H}_1\oplus\cdots\oplus\mathcal{H}_n)$;\\[3mm]
$(iii)$ $(a)$ $D_1\in\Phi_+(\mathcal{H}_1)$;\\
\hspace*{8mm}$(b)$ $\Big(\beta(D_j)=\infty$ for some $j\in\lbrace1,...,n\rbrace$ and $\alpha(D_s)<\infty$ for $2\leq s\leq j\Big)$ or $\Big(D_s\in\Phi_+(\mathcal{H}_s)$ for $2\leq s\leq n$ and $\sum\limits_{s=1}^n\alpha(D_s)\leq\sum\limits_{s=1}^n\beta(D_s)\Big)$.\\[3mm]
Then $(i) \Rightarrow (ii) \Rightarrow (iii)$.
\end{Teorema}
\begin{Primedba}
If $j=1$ in $(iii)(b)$, we simply omit condition "$\alpha(D_s)<\infty$ for $2\leq s\leq  j$" there.
\end{Primedba}

\textbf{Proof}: $(ii)\Rightarrow(iii)$

\hspace*{6mm}Assume that $T_n^d(A)$ is upper Weyl. Then $T_n^d(A)$ is upper Fredholm, hence $D_1\in\Phi_+(\mathcal{H}_1)$ (Lemma \ref{POMOCNALEMA}). Assume that $(iii)(b)$ fails. We have two possibilities. On the one hand, assume that for $2\leq s\leq n$ we have $\beta(D_s)<\infty$. It means (Theorem \ref{KATOZATVOREN}) that $\mathcal{R}(D_s)$ is closed for $1\leq s\leq n$. Again, we have two possibilities. Either there exists some $i\in\lbrace 2,...,n\rbrace$ with $\alpha(D_i)=\infty$, or we have $\sum\limits_{s=1}^n\alpha(D_s)>\sum\limits_{s=1}^n\beta(D_s)$.

First suppose $\alpha(D_i)=\infty$ for some $i\in\lbrace 2,...,n\rbrace.$ We use a method from \cite{WU3}. We know that for each $A\in\mathcal{B}_n$, operator matrix $T_n^d(A)$ as an operator from $\mathcal{N}(D_1)^\bot\oplus\mathcal{N}(D_1)\oplus\mathcal{N}(D_2)^\bot\oplus\mathcal{N}(D_2)\oplus\mathcal{N}(D_3)^\bot\oplus\mathcal{N}(D_3)\oplus\cdots\oplus\mathcal{N}(D_n)^\bot\oplus\mathcal{N}(D_n)$ into $\mathcal{R}(D_1)\oplus\mathcal{R}(D_1)^\bot\oplus\mathcal{R}(D_2)\oplus\mathcal{R}(D_2)^\bot\oplus\cdots\oplus\mathcal{R}(D_{n-1})\oplus\mathcal{R}(D_{n-1})^\bot\oplus\mathcal{R}(D_n)\oplus\mathcal{R}(D_n)^\bot$ admits the following block representation
\begin{equation}\label{MATRICA}
T_n^d(A)=\scalemath{0.85}
{\begin{bmatrix} 
    D_1^{(1)} & 0 & A_{12}^{(1)} & A_{12}^{(2)} & A_{13}^{(1)} & A_{13}^{(2)} & ... & A_{1n}^{(1)} & A_{1n}^{(2)}\\
    0 & 0 & A_{12}^{(3)} & A_{12}^{(4)} & A_{13}^{(3)} & A_{13}^{(4)} & ... & A_{1n}^{(3)} & A_{1n}^{(4)}\\
    0 & 0 & D_2^{(1)} & 0 & A_{23}^{(1)} & A_{23}^{(2)} & ... & A_{2n}^{(1)} & A_{2n}^{(2)}\\
    0 & 0 & 0 & 0 & A_{23}^{(3)} & A_{23}^{(4)} & ... & A_{2n}^{(3)} & A_{2n}^{(4)}\\
    0 & 0 & 0 & 0 & D_{3}^{(1)} & 0 & ... & A_{3n}^{(1)} & A_{3n}^{(2)}\\
    0 & 0 & 0 & 0 & 0 & 0 & ... & A_{3n}^{(3)} & A_{3n}^{(4)}\\
    \vdots & \vdots & \vdots & \vdots & \vdots & \vdots & \ddots & \vdots & \vdots\\
    0 & 0 & 0 & 0 & 0 & 0 & ... & A_{n-1,n}^{(1)} & A_{n-1,n}^{(2)}\\
    0 & 0 & 0 & 0 & 0 & 0 & ... & A_{n-1,n}^{(3)} & A_{n-1,n}^{(4)}\\
    0 & 0 & 0 & 0 & 0 & 0 & ... & D_n^{(1)} & 0\\
    0 & 0 & 0 & 0 & 0 & 0 & ... & 0 & 0\\
\end{bmatrix}}
\end{equation}
Obviously, $D_1^{(1)},\ D_2^{(1)},...,D_n^{(1)}$ from (\ref{MATRICA}) are invertible. Hence, there exist invertible operator matrices $U$ and $V$ so that 
\begin{equation}\label{MATRICA2}
UT_n^d(A)V=\scalemath{0.85}{\begin{bmatrix}
D_1^{(1)} & 0  & 0 & 0 & 0 & 0 & ... & 0 & 0\\
0 & 0 & 0 & B_{12}^{(4)} & 0 & B_{13}^{(4)} & ... & 0 & B_{1n}^{(4)}\\
0 & 0 & D_2^{(1)} & 0 & 0 & 0 & ... & 0 & 0\\
0 & 0 & 0 & 0 & 0 & B_{23}^{(4)} & ... & 0 & B_{2n}^{(4)}\\
0 & 0 & 0 & 0 & D_3^{(1)} & 0 & ... & 0 & 0\\
0 & 0 & 0 & 0 & 0 & 0 & ... & 0 & B_{3n}^{(4)}\\
\vdots & \vdots & \vdots & \vdots & \vdots & \vdots & \ddots & \vdots & \vdots\\
0 & 0 & 0 & 0 & 0 & 0 & ... & 0 & 0\\
0 & 0 & 0 & 0 & 0 & 0 & ... & 0 & B_{n-1,n}^{(4)}\\
0 & 0 & 0 & 0 & 0 & 0 & ... & D_n^{(1)} & 0\\
0 & 0 & 0 & 0 & 0 & 0 & ... & 0 & 0\\
\end{bmatrix}
}
\end{equation}
Next, it is clear that (\ref{MATRICA2}) is upper Weyl if and only if 
\begin{equation}\label{MATRICA3}
\begin{bmatrix}
0 & B_{12}^{(4)} & B_{13}^{(4)} & B_{14}^{(4)} & ... & B_{1n}^{(4)}\\
0 & 0          & B_{23}^{(4)} & B_{24}^{(4)} & ... & B_{2n}^{(4)}\\
0 & 0           & 0         & B_{34}^{(4)} & ... & B_{3n}^{(4)}\\
\vdots & \vdots   &   \vdots & \vdots & \ddots & \vdots\\
0 & 0          &   0        &    0      & ... & B_{n-1,n}^{(4)}\\
0 & 0 & 0 & 0 & ... &0
\end{bmatrix}
:
\begin{bmatrix}
\mathcal{N}(D_1)\\
\mathcal{N}(D_2)\\
\mathcal{N}(D_3)\\
\mathcal{N}(D_4)\\
\vdots\\
\mathcal{N}(D_n)
\end{bmatrix}
\rightarrow
\begin{bmatrix}
\mathcal{R}(D_1)^\bot\\
\mathcal{R}(D_2)^\bot\\
\mathcal{R}(D_3)^\bot\\
\vdots\\
\mathcal{R}(D_{n-1})^\bot\\
\mathcal{R}(D_{n})^\bot
\end{bmatrix}
\end{equation}
is upper Weyl. Since $\sum\limits_{s=1}^{i-1}\beta(D_s)<\infty$ and $\alpha(D_i)=\infty$, it follows that 
$$
\alpha\left(\begin{bmatrix}
B_{1i}^{(4)}\\
B_{2i}^{(4)}\\
B_{3i}^{(4)}\\
\vdots\\
B_{i-1,i}^{(4)}
\end{bmatrix}\right)=\infty,
$$
and hence operator defined in (\ref{MATRICA3}) is not upper Weyl for every $A\in\mathcal{B}_n$. This proves the desired.

Now assume $\alpha(D_s)<\infty$ for $2\leq s\leq n$. Then we have $\sum\limits_{s=1}^n\alpha(D_s)>\sum\limits_{s=1}^n\beta(D_s)$, and for each $A\in\mathcal{B}_n$, $T_n^d(A)$ has representation as (\ref{MATRICA}), and we still use (\ref{MATRICA2}) and (\ref{MATRICA3}). Since $D_s$, $1\leq s\leq n$ are upper Fredholm, we conclude that $T_n^d(A)$ is upper Weyl if and only if (\ref{MATRICA3}) is upper Weyl. From $\sum\limits_{s=1}^n\beta(D_s)<\sum\limits_{s=1}^n\alpha(D_s)$, we know (\ref{MATRICA3}) is not upper Weyl for every $A\in\mathcal{B}_n$.

On the other hand, assume that there is $j\in\lbrace 2
,...,n\rbrace$ with $\beta(D_j)=\infty$, and assume we have chosen the smallest such $j$. In that case $\beta(D_s)<\infty$ for $1\leq s\leq j-1$, hence $\mathcal{R}(D_s)$ is closed for $1\leq s\leq j-1$. Now, we easily conclude it is impossible that $\alpha(D_s)<\infty$ for $2\leq s\leq j-1$, otherwise $(iii)(b)$ would not fail. Therefore, $\alpha(D_j)=\infty$ for some $j\in\lbrace2,...,j-1\rbrace$ and be proceed with (\ref{MATRICA}), (\ref{MATRICA2}), (\ref{MATRICA3}) applied to $T_{j-1}^d(A)$.

$(i)\Rightarrow(ii)$

If $D_s\in\Phi_+(\mathcal{H}_s)$ for $2\leq s\leq n$ and $\sum\limits_{s=1}^n\alpha(D_s)\leq\sum\limits_{s=1}^n\beta(D_s)$, we trivially choose $A=(A_{ij}=)\mathbf{0}$. Assume that this is not the case. Otherwise, it holds  $\alpha(D_1)<\infty$, $\mathcal{R}(D_s)$ is closed for all $1\leq s\leq n$ and (\ref{PRVA}) holds. We find $A\in\mathcal{B}_n$ such that $\alpha(T_n^d(A))<\infty$ and $\mathcal{R}(T_n^d(A))$ is closed.  We choose $A=(A_{ij})_{1\leq i<j\leq n}$ so that $A_{ij}=0$ if $j-i\neq 1$, that is we place all nonzero operators of tuple $A$ on the superdiagonal. It remains to define $A_{ij}$ for $j-i=1$, $1\leq i<j\leq n$. First notice that $A_{i,i+1}:\mathcal{H}_{i+1}\rightarrow \mathcal{H}_i$. Since all of diagonal entries have closed ranges, we know that $\mathcal{H}_{i+1}=\mathcal{N}(D_{i+1})\oplus\mathcal{N}(D_{i+1})^\bot$, $\mathcal{H}_i=\mathcal{R}(D_i)^\bot\oplus\mathcal{R}(D_i)$, and from assumption (\ref{PRVA}) we get $\alpha(D_{i+1})\leq\beta(D_i)$. It follows that there is a left invertible operator $J_{i}:\mathcal{N}(D_{i+1})\rightarrow\mathcal{R}(D_i)^\bot$. We put $A_{i,i+1}=\begin{bmatrix}J_{i} & 0\\ 0 & 0\end{bmatrix}:\begin{bmatrix}\mathcal{N}(D_{i+1})\\ \mathcal{N}(D_{i+1})^\bot\end{bmatrix}=\mathcal{H}_{i+1}\rightarrow\mathcal{H}_i=\begin{bmatrix}\mathcal{R}(D_i)^\bot\\ \mathcal{R}(D_i)\end{bmatrix}$, and we implement this procedure for all $1\leq i\leq n-1$. Notice that $\mathcal{R}(D_i)$ is complemented to $\mathcal{R}(A_{i,i+1})$ for each $1\leq i\leq n-1.$

Now we have chosen our $A$, we show that $\mathcal{N}(T_n^d(A))\cong\mathcal{N}(D_1)$, implying $\alpha(T_n^d(A))=\alpha(D_1)<\infty$. Let us put $T_n^d(A)x=0$, where $x=x_1+\cdots+x_n\in \mathcal{H}_1\oplus\cdots\oplus \mathcal{H}_n$. The previous equality is then equivalent to the following system of equations
$$
\begin{bmatrix}D_1x_1+A_{12}x_2\\ D_2x_2+A_{23}x_3\\ \vdots\\ D_{n-1}x_{n-1}+A_{n-1,n}x_n\\ D_nx_n\end{bmatrix}=\begin{bmatrix}0\\ 0\\ \vdots\\ 0\\ 0\end{bmatrix}.
$$ 
The last equation gives $x_n\in\mathcal{N}(D_n)$. Since $\mathcal{R}(D_s)$ is complemented to $\mathcal{R}(A_{s,s+1})$ for all $1\leq s\leq n-1$, we have $D_sx_s=A_{s,s+1}x_{x+1}=0$ for all $1\leq s\leq n-1$. That is, $x_i\in\mathcal{N}(D_i)$ for every $1\leq i\leq n$, and $J_sx_{s+1}=0$ for every $1\leq s\leq n-1.$ Due to left invertibility of $J_s$ we get $x_s=0$ for $2\leq s\leq n$, which proves the claim. Therefore, $\alpha(T_n^d(A))=\alpha(D_1)<\infty$.

Secondly, we show that $\mathcal{R}(T_n^d(A))$ is closed. It is not hard to see that 
\begin{equation}\label{PRVA1}
\begin{aligned}
\mathcal{R}(T_n^d(A))=\mathcal{R}(D_1)\oplus\mathcal{R}(J_{1})\oplus\mathcal{R}(D_2)\oplus\mathcal{R}(J_{2})\oplus\cdots\oplus\mathcal{R}(D_{n-1})\oplus\\\mathcal{R}(J_{n-1})\oplus\mathcal{R}(D_n).
\end{aligned}
\end{equation} Furthermore, due to left invertibility of $J_i$'s, there exist closed subspaces $U_i$ of $\mathcal{R}(D_i)^\bot$ such that $\mathcal{R}(D_i)^\bot=\mathcal{R}(J_i)\oplus U_i$, $1\leq i\leq n-1$. It means that 
\begin{equation}\label{PRVA2}
\begin{aligned}
\mathcal{H}_1\oplus \mathcal{H}_2\oplus\cdots\oplus \mathcal{H}_n=\mathcal{R}(D_1)\oplus\mathcal{R}(J_{1})\oplus U_1\oplus\mathcal{R}(D_2)\oplus\mathcal{R}(J_{2})\oplus U_2\oplus\cdots\oplus\\\mathcal{R}(D_{n-1})\oplus\mathcal{R}(J_{n-1})\oplus U_{n-1}\oplus\mathcal{R}(D_n)\oplus\mathcal{R}(D_n)^\bot.
\end{aligned}
\end{equation} Comparing equalities (\ref{PRVA1}) and (\ref{PRVA2}), we conclude that $\mathcal{R}(T_n^d(A))$ is closed.

We have proved that $T_n^d(A)$ is upper Fredholm. Notice that $\beta(T_n^d(A))=\dim(U_1)+\dim(U_2)+\cdots+\dim(U_{n-1}) +\beta(D_n)$. Now, with respect to (\ref{PRVA}), either $\beta(D_n)=\infty$ or we can choose at least one $J_i$ such that its codimension is infinite, that is $\dim U_{i}=\infty$, $i\in\lbrace1,...,n-1\rbrace$. In either case we get $\beta(T_n^d(A))=\infty$ and it follows that $T_n^d(A)$ is upper Weyl. $\square$
\begin{Posledica}\label{POSLEDICA'}
Let $D_1\in\mathcal{B}(\mathcal{H}_1),\ D_2\in\mathcal{B}(\mathcal{H}_2),...,D_n\in\mathcal{B}(\mathcal{H}_n)$. Then
$$
\begin{aligned}
\sigma_{SF+}(D_1)\cup\Big(\bigcup\limits_{k=2}^{n+1}\Delta_k\Big)\subseteq\\\bigcap\limits_{A\in\mathcal{B}_n}\sigma_{aw}(T_n^d(A))\subseteq\\\sigma_{SF+}(D_1)\cup\Big(\bigcup_{k=2}^{n+1}\Delta_k'\Big)\cup\Big(\bigcup\limits_{k=2}^n\Delta_k''\Big),
\end{aligned}
$$
where
$$
\Delta_k:=\Big\lbrace\lambda\in\mathds{C}:\ \alpha(D_k-\lambda)=\infty\ and\  \sum\limits_{s=1}^{k-1}\beta(D_s-\lambda)<\infty\Big\rbrace,\ 2\leq k\leq n,
$$
$$
\Delta_{n+1}:=\Big\lbrace\lambda\in\mathds{C}:\ \sum\limits_{s=1}^n\beta(D_s-\lambda)<\sum\limits_{s=1}^{n}\alpha(D_s-\lambda)\Big\rbrace,
$$
$$
\Delta_k':=\lbrace\lambda\in\mathds{C}:\ \alpha(D_k-\lambda)>\beta(D_{k-1}-\lambda)\rbrace,\quad 2\leq k\leq n,
$$
$$
\Delta_{n+1}':=\Delta_{n+1},
$$
$$
\Delta_k'':=\Big\lbrace\lambda\in\mathds{C}:\ \mathcal{R}(D_k-\lambda)\ is\ not\ closed\Big\rbrace,\ 2\leq k\leq n.
$$
\end{Posledica}
\begin{Primedba}
Obviously, $\Delta_k\subseteq\Delta_k'$ for $2\leq k\leq n+1$.
\end{Primedba}
\begin{Teorema}\label{LEVIVEJLn=2'}
Let $D_1\in\mathcal{B}(\mathcal{H}_1), D_2\in\mathcal{B}(\mathcal{H}_2)$. Consider the following statements:\\[1mm]
$(i)$ $(a)$ $D_1\in\Phi_+(\mathcal{H}_1)$;\\
\hspace*{5.5mm}$(b)$ \Big($\alpha(D_2)\leq\beta(D_1)$, $\beta(D_1)+\beta(D_2)=\infty$ and $\mathcal{R}(D_2)$ is closed\Big) or $\Big(D_2\in\Phi_+(\mathcal{H}_2)\  and\ \alpha(D_1)+\alpha(D_2)\leq\beta(D_1)+\beta(D_2)\Big)$;\\[1mm]
$(ii)$ There exists $A\in\mathcal{B}_2$ such that $T_2^d(A)\in\Phi_+^-(\mathcal{H}_1\oplus\mathcal{H}_2)$;\\[1mm]
$(iii)$ $(a)$ $D_1\in\Phi_+(\mathcal{H}_1)$;\\
\hspace*{8.2mm}$(b)$ $\Big(\beta(D_1)=\infty$ or ($\beta(D_2)=\infty$ and $\alpha(D_1)<\infty)\Big)$ or \Big($D_2\in\Phi_+(\mathcal{H}_2)$ and $\alpha(D_1)+\alpha(D_2)\leq\beta(D_1)+\beta(D_2)$\Big).\\[1mm]
Then $(i)\Rightarrow(ii)\Rightarrow(iii)$.
\end{Teorema}
\begin{Posledica}\label{POSLEDICA2'}
Let $D_1\in\mathcal{B}(\mathcal{H}_1),\ D_2\in\mathcal{B}(\mathcal{H}_2)$. Then
$$
\sigma_{SF+}(D_1)\cup\Delta_2\cup\Delta_3\subseteq\\\bigcap\limits_{A\in\mathcal{B}_n}\sigma_{aw}(T_2^d(A))\subseteq\\\sigma_{SF+}(D_1)\cup\Delta_2'\cup\Delta_3\cup\Delta_2'',
$$
where
$$
\Delta_2:=\Big\lbrace\lambda\in\mathds{C}:\ \alpha(D_2-\lambda)=\infty\ and\ \beta(D_1-\lambda)<\infty\Big\rbrace,
$$
$$
\Delta_{3}:=\Big\lbrace\lambda\in\mathds{C}:\ \alpha(D_1-\lambda)+\alpha(D_2-\lambda)>\beta(D_1-\lambda)+\beta(D_2-\lambda)\Big\rbrace,
$$
$$
\Delta_2':=\lbrace\lambda\in\mathds{C}:\ \alpha(D_2-\lambda)\geq\beta(D_{1}-\lambda)\rbrace,
$$
$$
\Delta_2'':=\Big\lbrace\lambda\in\mathds{C}:\ \mathcal{R}(D_2-\lambda)\ is\ not\ closed\Big\rbrace.
$$
\end{Posledica}
\begin{Primedba}
Notice that $\Delta_2\subseteq\Delta_2'.$
\end{Primedba}

Statements concerning the lower Weyl spectrum of $T_n^d(A)$ we get by duality. 

\begin{Teorema}\label{DESNIVEJL'}
Let $D_1\in\mathcal{B}(\mathcal{H}_1),\ D_2\in\mathcal{B}(\mathcal{H}_2),...,D_n\in\mathcal{B}(\mathcal{H}_n)$. Consider the following statements:\\
$(i)$ $(a)$ $D_n\in\Phi_-(\mathcal{H}_n)$;\\
\hspace*{5.5mm}$(b)$ $\mathcal{R}(D_s)$ is closed for $1\leq s\leq n-1$ and 
\begin{equation}\label{DRUGA}
\begin{aligned}
\Big(\beta(D_s)\leq\alpha(D_{s+1})\quad for\ 1\leq s\leq n-1,\\
\sum_{s=1}^n\alpha(D_s)=\infty\Big)
\end{aligned}
\end{equation}
or \Big($D_s\in\Phi_-(\mathcal{H}_s)$ for $1\leq s\leq n-1$ and $\sum\limits_{s=1}^n\beta(D_s)\leq\sum\limits_{s=1}^n\alpha(D_s)$\Big);\\[3mm]
$(ii)$ There exists $A\in\mathcal{B}_n$ such that $T_n^d(A)\in\Phi_-^+(\mathcal{H}_1\oplus\cdots\oplus\mathcal{H}_n)$;\\[3mm]
$(iii)$ $(a)$ $D_n\in\Phi_-(\mathcal{H}_n)$;\\
\hspace*{8.2mm}$(b)$ \Big($\alpha(D_j)=\infty$ for some $j\in\lbrace2,...,n\rbrace$ and $\beta(D_s)<\infty$ for $j\leq s\leq n-1$\Big) or \Big($D_s\in\Phi_-(\mathcal{H}_s)$ for $1\leq s\leq n-1$ and $\sum\limits_{s=1}^n\beta(D_s)\leq\sum\limits_{s=1}^n\alpha(D_s)$\Big).\\[3mm]
Then $(i) \Rightarrow (ii) \Rightarrow (iii)$.
\end{Teorema}
\begin{Primedba}
If $j=n$ in $(iii)(b)$, we simply omit condition ''$\beta(D_s)<\infty$ for $j\leq s\leq n-1$''.
\end{Primedba}
\textbf{Proof.} The result immediately follows from Theorem \ref{LEVIVEJL'}, having in mind the statements of Lemma \ref{VEZA} and Lemma \ref{DODATNALEMA}. $\square$
\begin{Posledica}\label{POSLEDICA3'}
Let $D_1\in\mathcal{B}(\mathcal{H}_1),\ D_2\in\mathcal{B}(\mathcal{H}_2),...,D_n\in\mathcal{B}(\mathcal{H}_n)$. Then
$$
\begin{aligned}
\sigma_{SF-}(D_n)\cup\Big(\bigcup\limits_{k=1}^{n-1}\Delta_k\Big)\cup\Delta_{n+1}\subseteq\\\bigcap\limits_{A\in\mathcal{B}_n}\sigma_{sw}(T_n^d(A))\subseteq\\\sigma_{SF-}(D_n)\cup\Big(\bigcup\limits_{k=1}^{n-1}\Delta_k'\Big)\cup\Delta_{n+1}\cup\Big(\bigcup_{k=1}^{n-1}\Delta_k''\Big),
\end{aligned}
$$
where
$$
\Delta_k:=\Big\lbrace\lambda\in\mathds{C}:\ \beta(D_k-\lambda)=\infty\ and\ \sum\limits_{s=k+1}^{n}\alpha(D_s-\lambda)<\infty\Big\rbrace,\ 1\leq k\leq n-1,
$$
$$
\Delta_{n+1}:=\Big\lbrace\lambda\in\mathds{C}:\ \sum\limits_{s=1}^n\alpha(D_s-\lambda)<\sum\limits_{s=1}^n\beta(D_s-\lambda)\Big\rbrace,
$$
$$
\Delta_{k}':=\lbrace\lambda\in\mathds{C}:\ \beta(D_k-\lambda)>\alpha(D_{k+1}-\lambda)\rbrace,\quad 1\leq k\leq n-1,
$$
$$
\Delta_k'':=\Big\lbrace\lambda\in\mathds{C}:\ \mathcal{R}(D_k-\lambda)\ is\ not\ closed\Big\rbrace,\ 2\leq k\leq n-1.
$$
\end{Posledica}
\begin{Primedba}
Obviously, $\Delta_k\subseteq\Delta_k'$ for $1\leq k\leq n-1$.
\end{Primedba}
\begin{Teorema}\label{DESNIVEJLn=2'}
Let $D_1\in\mathcal{B}(\mathcal{H}_1),\ D_2\in\mathcal{B}(\mathcal{H}_2)$. Consider the following statements:\\
$(i)$ $(a)$ $D_2\in\Phi_-(\mathcal{H}_2)$;\\
\hspace*{5.5mm}$(b)$ \Big($\beta(D_1)\leq\alpha(D_{2})$, $\alpha(D_1)+\alpha(D_2)=\infty$ and $\mathcal{R}(D_1)$ is closed\Big) or \Big($D_1\in\Phi_-(\mathcal{H}_1)$ and $\beta(D_1)+\beta(D_2)\leq\alpha(D_1)+\alpha(D_2)$\Big);\\[3mm]
$(ii)$ There exists $A\in\mathcal{B}_2$ such that $T_2^d(A)\in\Phi_-^+(\mathcal{H}_1\oplus\mathcal{H}_2)$;\\[3mm]
$(iii)$ $(a)$ $D_2\in\Phi_-(\mathcal{H}_2)$;\\
\hspace*{8mm}$(b)$ \Big($\alpha(D_2)=\infty$ or $(\alpha(D_1)=\infty$ and $\beta(D_2)<\infty)\Big)$ or ($D_1\in\Phi_-(\mathcal{H}_1)$ and $\beta(D_1)+\beta(D_2)\leq\alpha(D_1)+\alpha(D_2)$).\\[3mm]
Then $(i) \Rightarrow (ii) \Rightarrow (iii)$.
\end{Teorema}

\begin{Posledica}\label{POSLEDICA4'}
Let $D_1\in\mathcal{B}(\mathcal{H}_1),\ D_2\in\mathcal{B}(\mathcal{H}_2)$. Then
$$
\sigma_{SF-}(D_2)\cup\Delta_1\cup\Delta_3\subseteq\\\bigcap\limits_{A\in\mathcal{B}_2}\sigma_{sw}(T_2^d(A))\subseteq\\\sigma_{SF-}(D_2)\cup\Delta_1'\cup\Delta_3\cup\Delta_1'',
$$
where
$$
\Delta_1:=\Big\lbrace\lambda\in\mathds{C}:\ \beta(D_1-\lambda)=\infty\ and\ \alpha(D_2-\lambda)<\infty\Big\rbrace,
$$
$$
\Delta_{3}:=\Big\lbrace\lambda\in\mathds{C}:\ \alpha(D_1-\lambda)+\alpha(D_2-\lambda)<\beta(D_1-\lambda)+\beta(D_2-\lambda)\Big\rbrace,
$$
$$
\Delta_1':=\lbrace\lambda\in\mathds{C}:\ \beta(D_1-\lambda)\geq\alpha(D_{2}-\lambda)\rbrace,
$$
$$
\Delta_1'':=\Big\lbrace\lambda\in\mathds{C}:\ \mathcal{R}(D_1-\lambda)\ is\ not\ closed\Big\rbrace.
$$
\end{Posledica}
\begin{Primedba}
Notice that $\Delta_1\subseteq\Delta_1'.$
\end{Primedba}

\section{Fredholmness of $T_n^d(A)$}
Now, we deal with characterizations of Fredholmness of $T_n^d(A)$. Since (upper, lower) Weyl operators form a subclass of (upper, lower) Fredholm operators, all theorems to follow will be reminiscent to theorems of Section 2. Therefore, corresponding proofs are a special case of the proofs already seen in the previous section, and so we omit them here.

We still assume $\mathcal{H}_1,...,\mathcal{H}_n$ to be arbitrary infinite dimensional Hilbert space. We generalize results of \cite{CAO},\cite{ZHANG} from $n=2$ to an arbitrary dimension of upper triangular operators, and we pose perturbation results of \cite{WU2} without assuming separability of underlying spaces. 

We start with a result which deals with the upper Fredholm spectrum of $T_n^d(A)$. 
\begin{Teorema}\label{LEVIFREDHOLM'}
Let $D_1\in\mathcal{B}(\mathcal{H}_1),\ D_2\in\mathcal{B}(\mathcal{H}_2),...,D_n\in\mathcal{B}(\mathcal{H}_n)$. Consider the following statements:\\
$(i)$ $(a)$ $D_1\in\Phi_+(\mathcal{H}_1)$;\\
\hspace*{5.5mm}$(b)$ $\mathcal{R}(D_s)$ is closed for $2\leq s\leq n$ and 
\begin{equation}\label{TRECA}
\begin{aligned}
\alpha(D_s)\leq\beta(D_{s-1})\quad for\ 2\leq s\leq n
\end{aligned}
\end{equation}
or $D_s\in\Phi_+(\mathcal{H}_s)$ for $2\leq s\leq n$;\\[3mm]
$(ii)$ There exists $A\in\mathcal{B}_n$ such that $T_n^d(A)\in\Phi_+(\mathcal{H}_1\oplus\cdots\oplus\mathcal{H}_n)$;\\[3mm]
$(iii)$ $(a)$ $D_1\in\Phi_+(\mathcal{H}_1)$;\\
\hspace*{8.2mm}$(b)$ \Big($\beta(D_j)=\infty$ for some $j\in\lbrace1,...,n-1\rbrace$ and $\alpha(D_s)<\infty$ for $2\leq s\leq j$\Big) or $D_s\in\Phi_+(\mathcal{H}_s)$ for $2\leq s\leq n$.\\[3mm]
Then $(i) \Rightarrow (ii) \Rightarrow (iii)$.
\end{Teorema}
\begin{Primedba}
If $j=1$ in $(iii)(b)$, we simply omit condition "$\alpha(D_s)<\infty$ for $2\leq s\leq  j$" there.
\end{Primedba}

\begin{Posledica}\label{POSLEDICA5'}
Let $D_1\in\mathcal{B}(\mathcal{H}_1),\ D_2\in\mathcal{B}(\mathcal{H}_2),...,D_n\in\mathcal{B}(\mathcal{H}_n)$. Then
$$
\begin{aligned}
\sigma_{SF+}(D_1)\cup\Big(\bigcup\limits_{k=2}^{n}\Delta_k\Big)\subseteq\\\bigcap\limits_{A\in\mathcal{B}_n}\sigma_{SF+}(T_n^d(A))\subseteq\\\sigma_{SF+}(D_1)\cup\Big(\bigcup_{k=2}^{n}(\Delta_k'\cap\Delta')\Big)\cup\Big(\bigcup\limits_{k=2}^n\Delta_k''\Big),
\end{aligned}
$$
where
$$
\Delta_k:=\Big\lbrace\lambda\in\mathds{C}:\ \alpha(D_k-\lambda)=\infty\ and\  \sum\limits_{s=1}^{k-1}\beta(D_s-\lambda)<\infty\Big\rbrace,\ 2\leq k\leq n,
$$
$$
\Delta_k':=\lbrace\lambda\in\mathds{C}:\ \alpha(D_k-\lambda)>\beta(D_{k-1}-\lambda)\rbrace,\quad 2\leq k\leq n,
$$
$$
\Delta':=\Big\lbrace\lambda\in\mathds{C}:\ \sum\limits_{s=2}^n\alpha(D_s-\lambda)=\infty\Big\rbrace,
$$
$$
\Delta_k'':=\Big\lbrace\lambda\in\mathds{C}:\ \mathcal{R}(D_k-\lambda)\ is\ not\ closed\Big\rbrace,\ 2\leq k\leq n.
$$
\end{Posledica}
\begin{Primedba}
Obviously, $\Delta_k\subseteq\Delta_k'\cap\Delta'$ for $2\leq k\leq n$.
\end{Primedba}
\begin{Teorema}\label{LEVIFREDHOLMn=2'}
Let $D_1\in\mathcal{B}(\mathcal{H}_1), D_2\in\mathcal{B}(\mathcal{H}_2)$. Consider the following statements:\\[1mm]
$(i)$ $(a)$ $D_1\in\Phi_+(\mathcal{H}_1)$;\\
\hspace*{5.5mm}$(b)$ \Big($\alpha(D_2)\leq\beta(D_1)$ and $\mathcal{R}(D_2)$ is closed\Big) or $D_2\in\Phi_+(\mathcal{H}_2)$;\\[1mm]
$(ii)$ There exists $A\in\mathcal{B}_2$ such that $T_2^d(A)\in\Phi_+(\mathcal{H}_1\oplus\mathcal{H}_2)$;\\[1mm]
$(iii)$ $(a)$ $D_1\in\Phi_+(\mathcal{H}_1)$;\\
\hspace*{8.2mm}$(b)$ $\beta(D_1)=\infty$ or $D_2\in\Phi_+(\mathcal{H}_2)$.\\[1mm]
Then $(i)\Rightarrow(ii)\Rightarrow(iii)$.
\end{Teorema}
\begin{Posledica}\label{POSLEDICA6'}
Let $D_1\in\mathcal{B}(\mathcal{H}_1),\ D_2\in\mathcal{B}(\mathcal{H}_2)$. Then
$$
\sigma_{SF+}(D_1)\cup\Delta_2\subseteq\\\bigcap\limits_{A\in\mathcal{B}_n}\sigma_{SF+}(T_2^d(A))\subseteq\\\sigma_{SF+}(D_1)\cup\Delta_2'\cup\Delta_2'',
$$
where
$$
\Delta_2:=\Big\lbrace\lambda\in\mathds{C}:\ \alpha(D_2-\lambda)=\infty\ and\ \beta(D_1-\lambda)<\infty\Big\rbrace,
$$
$$
\Delta_2':=\lbrace\lambda\in\mathds{C}:\ \alpha(D_2-\lambda)\geq\beta(D_{1}-\lambda)\rbrace,
$$
$$
\Delta_2'':=\Big\lbrace\lambda\in\mathds{C}:\ \mathcal{R}(D_2-\lambda)\ is\ not\ closed\Big\rbrace.
$$
\end{Posledica}
\begin{Primedba}
Notice that $\Delta_2\subseteq\Delta_2'.$
\end{Primedba}

Statements concerning the lower Fredholm spectrum of $T_n^d(A)$ we get by duality. 

\begin{Teorema}\label{DESNIFREDHOLM'}
Let $D_1\in\mathcal{B}(\mathcal{H}_1),\ D_2\in\mathcal{B}(\mathcal{H}_2),...,D_n\in\mathcal{B}(\mathcal{H}_n)$. Consider the following statements:\\
$(i)$ $(a)$ $D_n\in\Phi_-(\mathcal{H}_n)$;\\
\hspace*{5.5mm}$(b)$ $\mathcal{R}(D_s)$ is closed for $1\leq s\leq n-1$ and 
\begin{equation}\label{CETVRTA}
\begin{aligned}
\beta(D_s)\leq\alpha(D_{s+1})\quad for\ 1\leq s\leq n-1
\end{aligned}
\end{equation}
or $D_s\in\Phi_-(\mathcal{H}_s)$ for $1\leq s\leq n-1$;\\[3mm]
$(ii)$ There exists $A\in\mathcal{B}_n$ such that $T_n^d(A)\in\Phi_-(\mathcal{H}_1\oplus\cdots\oplus\mathcal{H}_n)$;\\[3mm]
$(iii)$ $(a)$ $D_n\in\Phi_-(\mathcal{H}_n)$;\\
\hspace*{8.2mm}$(b)$ \Big($\alpha(D_j)=\infty$ for some $j\in\lbrace2,...,n\rbrace$ and $\beta(D_s)<\infty$ for $j\leq s\leq n-1$\Big) or $D_s\in\Phi_-(\mathcal{H}_s)$ for $1\leq s\leq n-1$.\\[3mm]
Then $(i) \Rightarrow (ii) \Rightarrow (iii)$.
\end{Teorema}
\begin{Primedba}
If $j=n$ in $(iii)(b)$, we simply omit condition ''$\beta(D_s)<\infty$ for $j\leq s\leq n-1$'' there.
\end{Primedba}
\begin{Posledica}\label{POSLEDICA7'}
Let $D_1\in\mathcal{B}(\mathcal{H}_1),\ D_2\in\mathcal{B}(\mathcal{H}_2),...,D_n\in\mathcal{B}(\mathcal{H}_n)$. Then
$$
\begin{aligned}
\sigma_{SF-}(D_n)\cup\Big(\bigcup\limits_{k=1}^{n-1}\Delta_k\Big)\subseteq\\\bigcap\limits_{A\in\mathcal{B}_n}\sigma_{SF-}(T_n^d(A))\subseteq\\\sigma_{SF-}(D_n)\cup\Big(\bigcup\limits_{k=1}^{n-1}(\Delta_k'\cap\Delta'\Big)\cup\Big(\bigcup_{k=1}^{n-1}\Delta_k''\Big),
\end{aligned}
$$
where
$$
\Delta_k:=\Big\lbrace\lambda\in\mathds{C}:\ \beta(D_k-\lambda)=\infty\ and\ \sum\limits_{s=k+1}^{n}\alpha(D_s-\lambda)<\infty\Big\rbrace,\ 1\leq k\leq n-1,
$$
$$
\Delta_{k}':=\lbrace\lambda\in\mathds{C}:\ \beta(D_k-\lambda)>\alpha(D_{k+1}-\lambda)\rbrace,\quad 1\leq k\leq n-1,
$$
$$
\Delta':=\lbrace\lambda\in\mathds{C}:\ \sum_{s=1}^{n-1}\beta(D_s-\lambda)=\infty\rbrace,
$$
$$
\Delta_k'':=\Big\lbrace\lambda\in\mathds{C}:\ \mathcal{R}(D_k-\lambda)\ is\ not\ closed\Big\rbrace,\ 2\leq k\leq n-1.
$$
\end{Posledica}
\begin{Primedba}
$\Delta_k\subseteq\Delta_k'\cap\Delta'$ for $1\leq k\leq n-1$.
\end{Primedba}
\begin{Teorema}\label{DESNIFREDHOLMn=2'}
Let $D_1\in\mathcal{B}(\mathcal{H}_1),\ D_2\in\mathcal{B}(\mathcal{H}_2)$. Consider the following statements:\\
$(i)$ $(a)$ $D_2\in\Phi_-(\mathcal{H}_2)$;\\
\hspace*{5.5mm}$(b)$ \Big($\beta(D_1)\leq\alpha(D_{2})$ and $\mathcal{R}(D_1)$ is closed\Big) or $D_1\in\Phi_-(\mathcal{H}_1)$;\\[3mm]
$(ii)$ There exists $A\in\mathcal{B}_2$ such that $T_2^d(A)\in\Phi_-(\mathcal{H}_1\oplus\mathcal{H}_2)$;\\[3mm]
$(iii)$ $(a)$ $D_2\in\Phi_-(\mathcal{H}_2)$;\\
\hspace*{8.2mm}$(b)$ $\alpha(D_2)=\infty$ or $D_1\in\Phi_-(\mathcal{H}_1)$.\\[3mm]
Then $(i) \Rightarrow (ii) \Rightarrow (iii)$.
\end{Teorema}

\begin{Posledica}\label{POSLEDICA8'}
Let $D_1\in\mathcal{B}(\mathcal{H}_1),\ D_2\in\mathcal{B}(\mathcal{H}_2)$. Then
$$
\sigma_{SF-}(D_2)\cup\Delta_1\subseteq\\\bigcap\limits_{A\in\mathcal{B}_2}\sigma_{SF-}(T_2^d(A))\subseteq\\\sigma_{SF-}(D_2)\cup\Delta_1'\cup\Delta_1'',
$$
where
$$
\Delta_1:=\Big\lbrace\lambda\in\mathds{C}:\ \beta(D_1-\lambda)=\infty\ and\ \alpha(D_2-\lambda)<\infty\Big\rbrace,
$$
$$
\Delta_1':=\lbrace\lambda\in\mathds{C}:\ \beta(D_1-\lambda)\geq\alpha(D_{2}-\lambda)\rbrace,
$$
$$
\Delta_1'':=\Big\lbrace\lambda\in\mathds{C}:\ \mathcal{R}(D_1-\lambda)\ is\ not\ closed\Big\rbrace.
$$
\end{Posledica}
\begin{Primedba}
Notice that $\Delta_1\subseteq\Delta_1'.$
\end{Primedba}

Last topic is the class $\Phi(\mathcal{H}_1\oplus\cdots\oplus\mathcal{H}_n)$ and its corresponding essential spectrum.
\begin{Teorema}\label{FREDHOLM'}
Let $D_1\in\mathcal{B}(\mathcal{H}_1),\ D_2\in\mathcal{B}(\mathcal{H}_2),...,D_n\in\mathcal{B}(\mathcal{H}_n)$. Consider the following statements:\\
$(i)$  $(a)$ $D_1\in\Phi_+(\mathcal{H}_1)$ and  $D_n\in\Phi_-(\mathcal{H}_n)$;\\
\hspace*{5.5mm}(b) \Big($\mathcal{R}(D_s)$ is closed for $2\leq s\leq n-1$ and \big($\alpha(D_s)=\beta(D_{s-1})$ for $2\leq s\leq n$ or $\alpha(D_s)\leq\beta(D_{s-1})<\infty$ for $2\leq s\leq n$\big)\Big) or \Big($D_j\in\Phi_+(\mathcal{H}_j)$ for $2\leq j\leq n$ and $D_k\in\Phi_-(\mathcal{H}_k)$ for $1\leq k\leq n-1$\Big); \\[3mm]
$(ii)$ There exists $A\in\mathcal{B}_n$ such that $T_n^d(A)\in\Phi(\mathcal{H}_1\oplus\cdots\oplus\mathcal{H}_n)$;\\[3mm]
$(iii)$ $(a)$ $D_1\in\Phi_+(\mathcal{H}_1)$ and  $D_n\in\Phi_-(\mathcal{H}_n)$;\\
\hspace*{8.2mm}$(b)$ \Big($\beta(D_j)=\infty$ for some $j\in\lbrace1,...,n-1\rbrace$ and $\alpha(D_s)<\infty$ for $2\leq s\leq j$, $\alpha(D_k)=\infty$ for some $k\in\lbrace2,...,n\rbrace$, and $\beta(D_s)<\infty$ for $k\leq s\leq n-1$, $k>j$\Big) or \Big($D_j\in\Phi_+(\mathcal{H}_j)$ for $2\leq j\leq n$ and $D_k\in\Phi_-(\mathcal{H}_k)$ for $1\leq k\leq n-1$\Big).\\[3mm]
Then $(i) \Rightarrow (ii) \Rightarrow (iii)$.
\end{Teorema}
\begin{Primedba}
If $j=1$ and/or $k=n$ in $(iii)(b)$, condition that is ought to hold for $2\leq s\leq j$ and/or $k\leq s\leq n-1$ is omitted there.
\end{Primedba}
\textbf{Proof. } $(ii) \Rightarrow (iii)$

Let $T_n^d(A)$ be Fredholm for some $A\in\mathcal{B}_n$. Then $T_n^d(A)$ is both left and lower Fredholm, and so by employing Theorems \ref{LEVIFREDHOLM'} and \ref{DESNIFREDHOLM'} we easily get the desired.

$(i)\Rightarrow(ii)$

If $D_j\in\Phi_+(\mathcal{H}_j)$ for $2\leq j\leq n$ and $D_k\in\Phi_-(\mathcal{H}_k)$ for $1\leq k\leq n-1$ choose trivially $A=\mathbf{0}$. Otherwise, this part follows the argument as seen in the proof of Theorem \ref{LEVIVEJL'}. Namely, assumptions of $(i)(b)$ ensure the existence of left invertible $J_i$'s, and so we choose $A=(A_{ij})$ as shown there. We shall again have $\alpha(T_n^d(A))=\alpha(D_1)<\infty$, and due to our assumptions we can choose all $U_i$'s to be finite dimensional. Therefore, $\beta(T_n^d(A))=\dim U_1+\cdots+\dim U_{n-1}+\beta(D_n)<\infty$, having in mind that $D_n\in\Phi_-(\mathcal{H}_n)$. $\square$

\begin{Posledica}\label{POSLEDICA9'}
Let $D_1\in\mathcal{B}(\mathcal{H}_1),\ D_2\in\mathcal{B}(\mathcal{H}_2),...,D_n\in\mathcal{B}(\mathcal{H}_n)$. Then
$$
\begin{aligned}
\sigma_{SF+}(D_1)\cup\sigma_{SF-}(D_n)\cup\Big(\bigcup\limits_{k=2}^{n-1}\Delta_k\Big)\cup\Delta_n\subseteq\\\bigcap_{A\in\mathcal{B}_n}\sigma_e(T_n^d(A))\subseteq\\\sigma_{SF+}(D_1)\cup\sigma_{SF-}(D_n)\cup\Big(\Big(\bigcup\limits_{k=2}^{n}\delta_k'\Big)\cap\Big(\bigcup\limits_{k=2}^{n}\Delta_k'\Big)\Big)\cup\Big(\bigcup\limits_{k=2}^{n-1}\Delta_k''\Big),
\end{aligned}
$$
where
$$
\begin{aligned}
\Delta_k=\Big\lbrace\lambda\in\mathds{C}:\ \alpha(D_k-\lambda)=\infty\ and\ \sum_{s=1}^{k-1}\beta(D_s-\lambda)<\infty\Big\rbrace\cup\\
\Big\lbrace\lambda\in\mathds{C}:\ \beta(D_k-\lambda)=\infty\ and\ \sum_{s=k+1}^n\alpha(D_s-\lambda)<\infty\Big\rbrace,\quad 2\leq k\leq n-1,
\end{aligned}
$$
$$
\begin{aligned}
\Delta_n=\Big\lbrace\lambda\in\mathds{C}:\ \alpha(D_n-\lambda)=\infty\ and\ \sum_{s=1}^{n-1}\beta(D_s-\lambda)<\infty\Big\rbrace\cup\\\Big\lbrace\lambda\in\mathds{C}:\ \beta(D_1-\lambda)=\infty\ and\ \sum_{s=2}^n\alpha(D_s-\lambda)<\infty\Big\rbrace,
\end{aligned}
$$
$$
\delta_k':=\lbrace\lambda\in\mathds{C}:\ \beta(D_{k-1}-\lambda)=\infty\ or\ \alpha(D_k-\lambda)>\beta(D_{k-1}-\lambda)\rbrace,\ 2\leq k\leq n,
$$
$$
\Delta_k':=\lbrace\lambda\in\mathds{C}:\ \alpha(D_s-\lambda)\neq\beta(D_{s-1}-\lambda)\rbrace,\ 2\leq k\leq n,
$$
$$
\Delta_k'':=\lbrace\lambda\in\mathds{C}:\ \mathcal{R}(D_s-\lambda)\ is\ not\ closed\rbrace,\ 2\leq k\leq n-1.
$$
\end{Posledica}

\begin{Primedba}
Obviously, $\Delta_k\subseteq\Delta_k'\cap\delta_k'$ for each $2\leq k\leq n$.
\end{Primedba}

\begin{Teorema}\label{FREDHOLMn=2'}
Let $D_1\in\mathcal{B}(\mathcal{H}_1),\ D_2\in\mathcal{B}(\mathcal{H}_2)$. Consider the following statements:\\
$(i)$  $(a)$ $D_1\in\Phi_+(\mathcal{H}_1)$ and  $D_2\in\Phi_-(\mathcal{H}_2)$;\\
\hspace*{5.5mm}(b) \Big($\alpha(D_2)=\beta(D_1)$ or $\alpha(D_2)\leq\beta(D_{1})<\infty$\Big) or \Big($D_2\in\Phi_+(\mathcal{H}_2)$ and $D_1\in\Phi_-(\mathcal{H}_1)$\Big). \\[3mm]
$(ii)$ There exists $A\in\mathcal{B}_2$ such that $T_2^d(A)\in\Phi(\mathcal{H}_1\oplus\mathcal{H}_2)$.\\[3mm]
$(iii)$ $(a)$ $D_1\in\Phi_+(\mathcal{H}_1)$ and  $D_2\in\Phi_-(\mathcal{H}_2)$;\\
\hspace*{8.2mm}$(b)$ \Big($\alpha(D_2)=\beta(D_1)=\infty$\Big) or \Big($D_2\in\Phi_+(\mathcal{H}_2)$ and $D_1\in\Phi_-(\mathcal{H}_1)$\Big).\\[3mm]
Then $(i) \Rightarrow (ii) \Rightarrow (iii)$.
\end{Teorema}

\begin{Posledica}\label{POSLEDICA10'}
Let $D_1\in\mathcal{B}(\mathcal{H}_1),\ D_2\in\mathcal{B}(\mathcal{H}_2)$. Then
$$
\sigma_{SF+}(D_1)\cup\sigma_{SF-}(D_2)\cup\Delta\subseteq\\\bigcap_{A\in\mathcal{B}_2}\sigma_e(T_2^d(A))\subseteq\\\sigma_{SF+}(D_1)\cup\sigma_{SF-}(D_2)\cup\Delta',
$$
where
$$
\begin{aligned}
\Delta=\lbrace\lambda\in\mathds{C}:\ \alpha(D_2-\lambda)=\infty\ and\ \beta(D_1-\lambda)<\infty\rbrace\cup\\\lbrace\lambda\in\mathds{C}:\ \beta(D_1-\lambda)=\infty\ and\ \alpha(D_2-\lambda)<\infty\rbrace,
\end{aligned}
$$
$$
\begin{aligned}
\Delta'=\lbrace\lambda\in\mathds{C}:\ \alpha(D_2-\lambda)\neq\beta(D_1-\lambda)\rbrace\cap\\\lbrace\lambda\in\mathds{C}:\ \beta(D_1-\lambda)=\infty\ or\ \alpha(D_2-\lambda)>\beta(D_1-\lambda)\rbrace.
\end{aligned}
$$
\end{Posledica}
\begin{Primedba}
Notice that $\Delta\subseteq\Delta'.$
\end{Primedba}

\end{document}